\documentclass[reqno, 12 pt]{amsart}
\usepackage{xcolor,fullpage,hyperref,enumitem}
\usepackage{tikz,amsmath,amssymb,color,mathtools}
\usepackage{graphicx}
\usepackage[margin=1in]{geometry}
\usepackage{amsthm}
\usetikzlibrary{calc}
\usepackage{float}
\usepackage{amsrefs}
\usetikzlibrary{arrows.meta, positioning}
\usepackage{subcaption}

\theoremstyle{plain}
\newtheorem{theorem}{Theorem}

\newtheorem{proposition}[theorem]{Proposition}
\theoremstyle{definition}
\newtheorem{definition}[theorem]{Definition}

\theoremstyle{remark}
\newtheorem{remark}[theorem]{Remark}

\newcommand{\yellowcell}[2]{%
  \fill[yellow!70] (#1,#2) rectangle ++(1,1);
  \draw[thick]     (#1,#2) rectangle ++(1,1);}
\newcommand{\greencell}[2]{%
  \fill[green!70] (#1,#2) rectangle ++(1,1);
  \draw[thick]     (#1,#2) rectangle ++(1,1);}
\newcommand{\redcell}[2]{%
  \fill[red!70] (#1,#2) rectangle ++(1,1);
  \draw[thick]     (#1,#2) rectangle ++(1,1);}
\newcommand{\bluecell}[2]{%
  \fill[blue!70] (#1,#2) rectangle ++(1,1);
  \draw[thick]     (#1,#2) rectangle ++(1,1);}

\def\GFt{GFt}

\newcommand{\seqnum}[1]{\href{https://oeis.org/#1}{\rm \underline{#1}}}
\usepackage{thmtools} 
\usepackage{thm-restate}

\title{Automatic Enumeration of Tilings by Polyominoes}

\author[Martinez]{Lucy Martinez}
\address[L.~Martinez]{Department of Mathematics, Rutgers University, Piscataway, NJ 08854}
\email{\textcolor{blue}{\href{mailto:lucy.martinez@rutgers.edu}{lucy.martinez@rutgers.edu}}}

\begin{document}

\begin{abstract}
We revisit Zeilberger's computational framework for enumerating polyomino tilings and demonstrate its use in counting tilings of $k\times n$ boards by $L$-tetrominoes, allowing rotations. Recently, B\v{e}lohoubek and Slav\'ik gave a generating function for tilings of $2n\times 4$ boards by $L$-tetrominoes. Motivated by their work, we first verify the method by automatically recovering their generating function. We then apply the method to widths $5$ through $9$ and obtain explicit rational generating functions and new integer sequences that are not currently listed in the OEIS.
\end{abstract}

\maketitle

\section{Introduction}

A \emph{polyomino} is a connected geometric shape formed by joining a finite number of unit cells edge-to-edge. A polyomino of size $n$ (formed by exactly $n$ cells) is called an \emph{$n$-omino}. The word polyomino was invented by Solomon Golomb in 1953 when he gave a talk to the Harvard Mathematics Club~\cite{Golomb94}. In 1965, Golomb published a book that included a collection of his previous results on polyomino problems, and in the 1994 revised edition, Golomb included subsequent developments since the first edition. Polyominoes have been used in popular puzzles, but they also appear in statistical physics through polymer models and connected clusters~\cite{Polymers79,Read1962}. In combinatorics, one of the problems that is still open is enumerating polyominoes of a given size. Only special classes of polyominoes have a known formula.

A classical problem involving polyominoes is tiling a board with a specific set of polyominoes. Dominoes arise in several tiling problems, including the Aztec diamond problem, in which the number of domino tilings of the Aztec diamond is given by a power of $2$~\cite{EKLP}, and the mutilated chessboard problem in which two opposite corners are removed, making it impossible to tile such a board with $31$ dominoes~\cite{Gardner57}. In general, obtaining exact formulas for the number of tilings of a $k\times n$ board is a difficult problem. But in 1961, Temperley and Fisher~\cite{Dimer}, and Kasteleyn~\cite{Kasteleyn1987} independently solved the celebrated dimer problem~\cite{Dimer}.
They gave an explicit formula for the number of tilings of a $k\times n$ board by $2\times 1$ dominoes.
For $m\geq 3$, a formula that enumerates the number of tilings of $k\times n$ boards by $m\times 1$ tiles with arbitrary integers $k,n$ is not known. Klarner and Pollack~\cite{KlarnerPollack} and Stanley~\cite{Stanley85} subsequently studied domino tilings of rectangular boards with fixed width and developed recurrences and generating functions for their enumeration.

\begin{figure}[H]
    % First subfigure
    \begin{subfigure}{0.15\textwidth}
    \centering
        \begin{tikzpicture}[x=0.6cm,y=0.6cm]
           \greencell{0}{0}
           \greencell{0}{1}
        \end{tikzpicture}
    \caption*{Free domino}
    \end{subfigure}
    % Second subfigure
    \begin{subfigure}{0.25\textwidth}
    \centering
        \begin{tikzpicture}[x=0.6cm,y=0.6cm]
            \yellowcell{0}{0}
            \yellowcell{0}{1}
            \yellowcell{0}{2}
            \yellowcell{2}{1}
            \yellowcell{2}{0}
            \yellowcell{3}{0}
        \end{tikzpicture}
    \caption*{Free trominoes}
    \end{subfigure}
    % Third subfigure
    \begin{subfigure}{0.5\textwidth}
    \centering
        \begin{tikzpicture}[x=0.6cm,y=0.6cm]
            %I shape:
            \redcell{0}{0}
            \redcell{0}{1}
            \redcell{0}{2}
            \redcell{0}{3}
            %L-shape:
            \redcell{2}{2}
            \redcell{2}{1}
            \redcell{2}{0}\redcell{3}{0}
            %S- shape
            \redcell{5}{2}
            \redcell{5}{1}\redcell{6}{1}
            \redcell{6}{0}
            %T- shape
            \redcell{8}{2}
            \redcell{8}{1}\redcell{9}{1}
            \redcell{8}{0}
            %O shape
            \redcell{11}{1}\redcell{12}{1}
            \redcell{11}{0}\redcell{12}{0}
        \end{tikzpicture}
    \caption*{Free tetrominoes}
    \end{subfigure}
    \caption{Free polyominoes with $2,3$ and $4$ cells.}
    \label{fig:freepolys1}
\end{figure}
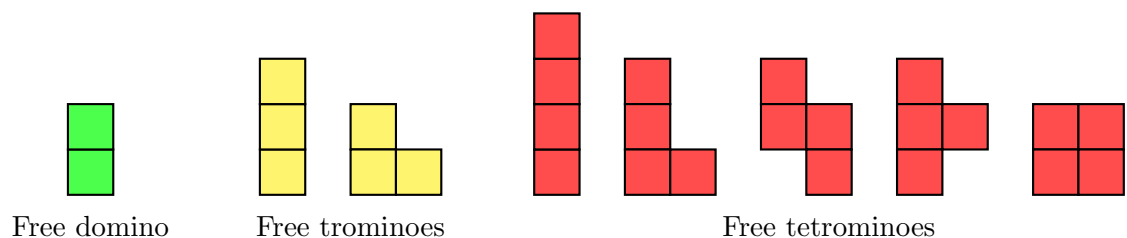

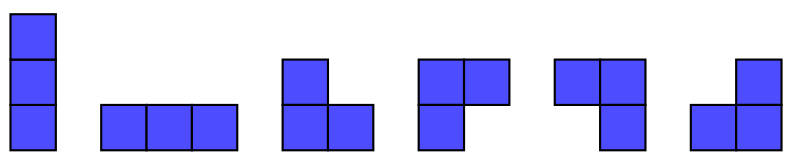
\begin{figure}[H]
    \begin{tikzpicture}[x=0.6cm,y=0.6cm]
        \bluecell{0}{0}
        \bluecell{0}{1}
        \bluecell{0}{2}
        %second shape
        \bluecell{2}{0}
        \bluecell{3}{0}
        \bluecell{4}{0}
        %L-shape:
        \bluecell{6}{1}
        \bluecell{6}{0}\bluecell{7}{0}
        %fourth
        \bluecell{9}{1}\bluecell{10}{1}
        \bluecell{9}{0}
        %fifth
        \bluecell{12}{1}\bluecell{13}{1}
        \bluecell{13}{0}
        %sixth
        \bluecell{16}{1}
        \bluecell{15}{0}\bluecell{16}{0}
    \end{tikzpicture}
    \caption{The six fixed trominoes, which allow only translation.}
    \label{fig:fixedpolys1}
\end{figure}

Before discussing other polyomino tiling problems, we recall the standard equivalence classes used for polyominoes. The two main classes of polyominoes are free polyominoes and fixed polyominoes. The main transformations considered for polyominoes are translation, rotation, and reflection. Free polyominoes are polyominoes considered up to translation, rotation, and reflection; that is, two free polyominoes are considered equivalent if one can be obtained from the other by a combination of these transformations. In contrast, fixed polyominoes are considered distinct under rotations and reflections. In Figure~\ref{fig:freepolys1}, the families of free polyominoes with $2, 3$, and $4$ cells are shown. In Figure~\ref{fig:fixedpolys1}, all six fixed polyominoes with $3$ cells are shown. Moreover, some polyominoes have names. For example, in Figure~\ref{fig:freepolys1}, from left to right, they are the straight tetromino, the $L$-tetromino (or right tetromino), the skew tetromino, the $T$-tetromino, and the square tetromino.

Throughout this paper, unless stated otherwise, rotations of the tiles are allowed, but reflections are not.
In Golomb's book~\cite{Golomb94}, one of the open problems was to determine the number of tilings of a $4\times n$ board by $L$-trominoes, allowing rotations only.
In 1982, Read~\cite{Read1982} introduced a unified method for enumerating tilings of narrow rectangular arrays using a frontier representation of the partially tiled board. Although Read's work focused on squares and dominoes, the method applies more generally and provides an important approach used in algorithmic enumeration of polyomino tilings. In 1999, Moore solved Golomb's problem by providing generating functions for tilings of $4 \times n$ and $5 \times n$ boards by $L$-trominoes, and tilings of $4 \times n$ boards by $L$- and $T$-tetrominoes, under the same convention~\cite{Moore1999}.
In 2006, Zeilberger developed a Maple (a computer algebra system) package that automatically constructs generating functions for tiling a $k\times n$ board with fixed width using an arbitrary finite set of tiles~\cite{DZ}. The algorithm generates finitely many auxiliary ``jagged'' boards, or states, obtained by placing tiles along the frontier, and constructs a system of linear equations for their generating functions. This finite state construction is closely related to the transfer-matrix method (see Stanley~\cite[Section 4.7]{StanleyVol1}). Related finite state approaches have also been used to enumerate tilings of rectangular boards. For example, Mathar~\cite{Mathar} used the transfer-matrix to obtain generating functions for tiling rectangular boards using rectangular tiles.

In a recent paper in the American Mathematical Monthly, B\v{e}lohoubek and Slav\'ik give explicit formulas for the number of tilings of a $2n\times 4$ board by $L$-tetrominoes, allowing rotations~\cite{Ltetromino}. They also show that such tilings are equinumerous with a family of two-color integer compositions. Their results motivate the present paper.
This paper revisits Zeilberger's automated framework and demonstrates its effectiveness for enumerating polyomino tilings. In Section~\ref{subsec:Ltetro}, we show that the generating function provided by B\v{e}lohoubek and Slav\'ik can be recovered automatically from Zeilberger's construction. We then apply this framework to obtain generating functions for the number of tilings of $k\times n$ boards by $L$-tetrominoes, where $k\in \{5,6,7,8,9\}$. These sequences appear to be new and suitable for inclusion in the OEIS~\cite{OEIS}. 
This computational approach is not limited to this specific setting, but it also produces explicit generating functions for larger widths (at least up to $k=16$) with an arbitrary set of tiles. Although the work of B\v{e}lohoubek and Slav\'ik provides combinatorial insight, our focus is on automatic enumeration via computational methods and generating functions.
The Maple code used to produce the computations in this paper is publicly available on the GitHub repository~\cite{github}.

\section{Program Framework for tiling a board}\label{sec:program}

\begin{figure}[H]
   \begin{tikzpicture}[x=0.6cm,y=0.6cm]
        \redcell{0}{0}
        \redcell{0}{1}
        \draw[->, thick] (0,0) -- (2.5,0) node[right] {$x$};
        \draw[->, thick] (0,0) -- (0,3.5) node[above] {$y$};
        \foreach \x in {0,1,2}
        {
        \draw (\x,0.1) -- (\x,-0.1);
        \node[below] at (\x,0) {\small $\x$};
        }
        \foreach \y in {1,2,3}
        {
        \draw (0.1,\y) -- (-0.1,\y);
        \node[left] at (0,\y) {\small $\y$};
        }
    \end{tikzpicture}
    \caption{The set of lattice points $\{(0,0), (0,1)\}$ of the $2\times 1$ domino.}
    \label{fig:domino}
\end{figure}
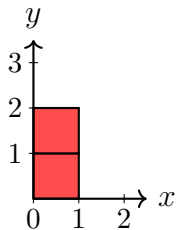

Zeilberger developed a Maple program that produces the generating function for the number of tilings of a $k\times n$ board by a given set of tiles, where $k$ is fixed and $n$ is arbitrary~\cite{DZ}. We give an overview of the program and include new sequences not yet in the OEIS.
We represent each polyomino with a normalized set of lattice points. We use the coordinate system with the origin at the bottom left corner of the $k\times n$ board, such that one side of length $k$ of the board is along the $y$-axis. Each cell of a polyomino is represented by a pair $(x,y)$ corresponding to the coordinates of its bottom left corner. For instance, in Figure~\ref{fig:domino}, the free domino consists of two cells represented by the set $\{(0,0), (0,1)\}$.

\begin{definition}\label{def:generatingfunc}
Let $T$ be a collection of polyominoes such that $T=\{P_1,P_2,\ldots,P_m\}$ where each $P_i$ is a normalized set of lattice points. Define $\GFt(k,T,t)$ to be the generating function in the variable $t$ whose coefficient of $t^n$ is the number of tilings of the $k\times n$ board using tiles from the set $T$. We refer to $T$ as the \emph{set of tiles}.

For a positive integer $c$, we define $F(k,T,t,c)$ to be the generating function whose coefficient of $t^n$ is the number of tilings of a $k\times cn$ board using tiles from the set $T$.
\end{definition}

A classical example of tiling is the number of ways to tile a $2\times n$ board by the rotations of a domino. Using our notation, if $T=\{\{(0,0),(0,1)\},\{(0,0),(1,0)\}\}$ then $\GFt(2,T,t)=\frac{1}{1-t-t^2}$, which is the generating function of the Fibonacci numbers (OEIS \seqnum{A000045}). 
This simple example illustrates the type of output produced by Zeilberger's program for an arbitrary set of tiles and fixed board width. More generally, Zeilberger's algorithm constructs a finite system of linear equations for the generating functions associated with the original rectangular board and its ``auxiliary jagged'' boards. Solving this system yields exact rational generating functions. Since every rational generating function determines a linear recurrence for its coefficients, the recurrences follow immediately from the denominators. Indeed, if
\[
G(t)=\sum_{n\ge0}a_nt^n=\frac{P(t)}{Q(t)},
\]
where 
\[Q(t)=1+c_1t+\cdots+c_dt^d,\]
then $Q(t)G(t)=P(t)$. Comparing coefficients of $t^n$ for $n>\deg(P)$ yields
\[a_n=-c_1a_{n-1}-\cdots-c_da_{n-d}.\]

Thus, Zeilberger's program automatically computes the generating function and the corresponding recurrence for a fixed board width and a set of tiles.

We now illustrate this framework on the five free tetrominoes before turning to the $L$-tetromino, which is the primary focus of this paper.
Among them, the square and skew tetrominoes contribute little in our setting. There is only one way to tile a rectangular board by the square tetromino, while the skew tetromino cannot tile a rectangular board without leaving holes. 
We therefore begin with the straight tetromino, whose tiling sequences are largely known and provide a useful benchmark for Zeilberger's program. The first known sequences are summarized in Table~\ref{tab:straight-tetro-seqs}.

\begin{remark}\label{rem:convention}
We index each sequence using the smallest positive scaling factor that avoids the zero terms. That is, for each fixed $k$, we choose the smallest positive integer $c$ such that the number of tilings of a $k\times cn$ board is nonzero for every $n\geq 0$. The corresponding sequence is generated by $F(k,T,t,c)$, as described in Definition~\ref{def:generatingfunc}.
\end{remark}

\begin{table}[H]
    \centering
    \begin{tabular}{|c|c|c|}
     \hline 
       Board dimension  & First few terms of the sequence & OEIS \\\hline 
        $4\times n$ & $1, 1, 1, 1, 2, 3, 4, 5, 7, 10, 14, 19, 26, 36, 50$ &\seqnum{A003269} \\ \hline 
        $5\times 4n$ & $1, 3, 15, 75, 371, 1833, 9057, 44753, 221137$ &\seqnum{A236579}\\ \hline 
        $6\times 4n$ & $1, 4, 25, 154, 943, 5773, 35344, 216388, 1324801$ &\seqnum{A236580} \\ \hline 
        $7 \times 4n$ & $1,5,37,269,1949,14121,102313,741305,5371097$ &  \seqnum{A236581}\\ \hline 
        $8\times n$ & $1, 1, 1, 1, 7, 15, 25, 37, 100, 229, 454, 811, 1732, 3777$ &\seqnum{A236582} \\ \hline
    \end{tabular}
    \caption{Sequences of tiling boards with rotations of the straight tetromino.}
    \label{tab:straight-tetro-seqs}
\end{table}

\subsection{Tiling boards with $T$-tetrominoes}\label{subsec:Ttetro}

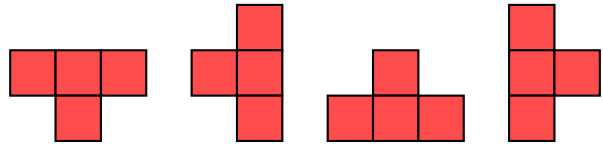
\begin{figure}[H]
\begin{tikzpicture}[x=0.6cm,y=0.6cm]
%T-shape
\redcell{0}{1}\redcell{1}{1}\redcell{2}{1}
\redcell{1}{0}
%second
\redcell{5}{2}
\redcell{4}{1}\redcell{5}{1}
\redcell{5}{0}
%third
\redcell{8}{1}
\redcell{7}{0}\redcell{8}{0}\redcell{9}{0} 
%fourth
\redcell{11}{2}
\redcell{11}{1}\redcell{12}{1}
\redcell{11}{0}
\end{tikzpicture}
\caption{The set of all rotations of the $T$-tetromino.}
\label{fig:Ttetro-rot}
\end{figure}

We consider the rotations of the $T$-tetromino, represented by \[P=\{(0, 1), (1, 0), (1, 1), (2, 1)\}.\] In Figure~\ref{fig:Ttetro-rot}, we display all the rotations of the $T$-tetromino.

In 1999, Moore provided a formula for the number of tilings of a $4\times 4 n$ board by $T$-tetrominoes~\cite{Moore1999}. The first few corresponding terms are (OEIS~\seqnum{A008776})
\[1,2, 6, 18, 54, 162, 486, 1458, 4374, 13122, 39366, \ldots\]
and the corresponding generating function is
\[F(4,T,t,4)=\frac{t-1}{3t-1},\]
where 
\begin{align*}
T=\{ &\{(0, 1), (1, 0), (1, 1), (2,1)\}, \{(0, 1), (1, 0), (1, 1), (1, 2)\}, \\
& \{(0, 0), (1, 0), (1, 1), (2, 0)\},  \{(0, 0), (0, 1), (0, 2), (1, 1)\} \}.
\end{align*}

If $k=5$, one cannot tile such a board with $T$-tetrominoes. Walkup showed that a board can be tiled with $T$-tetrominoes if and only if its length and width are both multiples of $4$~\cite{Walkup}. Consequently, we consider the sequences indexed by boards of dimensions $k\times 4n$, where $k\in \{8,12,16\}$.
If $\GFt(8, T,t)$ is the generating function for the number of tilings of an $8\times n$ board with $T$-tetrominoes, then Zeilberger's program outputs
\[\GFt(8, T, t)= \frac{15 t^{8}-10 t^4 +1}{27 t^{8}-16 t^4 +1},\]
which we rewrite as $F(8,T,t,4)$, the generating function for the $8\times 4n$ boards,
\[F(8, T,t,4)= \frac{15 t^{2}-10 t +1}{27 t^{2}-16 t +1}.\]

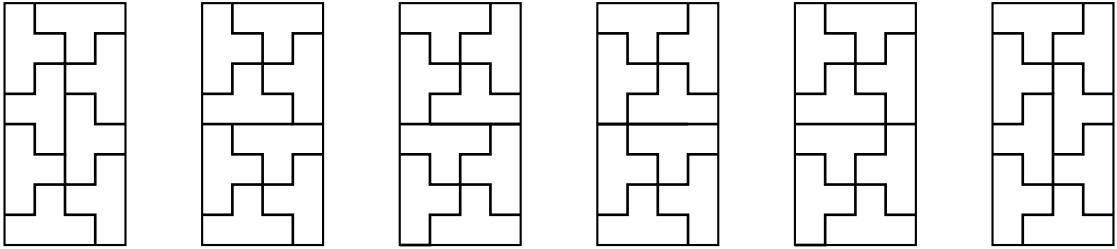
\begin{figure}[H]
    % First subfigure
    \begin{subfigure}{0.15\textwidth}
    \centering
    \begin{tikzpicture}[scale=0.4]
    % board
    \draw[thick] (0,0) rectangle (4,8);

    % Tile 1 (top left)
    \draw[line width=1pt]
    (0,5)--(1,5)--(1,6)--(2,6)--(2,7)--(1,7)--(1,8);

    \draw[line width=1pt]
     (2,6)--(3,6)--(3,7)--(4,7);

    \draw[line width=1pt]
    (0,4)--(1,4)--(1,3)--(2,3)--(2,6);

    \draw[line width=1pt]
     (2,5)--(3,5)--(3,4)--(4,4);

    \draw[line width=1pt]
    (2,3)--(2,2)--(3,2)--(3,3)--(4,3);
  
    \draw[line width=1pt]
    (0,1)--(1,1)--(1,2)--(2,2)--(2,1)--(3,1)--(3,0);
    \end{tikzpicture}
    \end{subfigure}
    % Second subfigure
    \begin{subfigure}{0.15\textwidth}
    \centering
     \begin{tikzpicture}[scale=0.4]
    % board
    \draw[thick] (0,0) rectangle (4,8);

    % Tile 1 (top left)
    \draw[line width=1pt]
    (0,5)--(1,5)--(1,6)--(2,6)--(2,7)--(1,7)--(1,8);

    \draw[line width=1pt]
     (2,6)--(3,6)--(3,7)--(4,7);

    \draw[line width=1pt]
     (2,6)--(2,5)--(3,5)--(3,4)--(4,4); 
     %%%%%%
     \draw[line width=1pt]
     (0,4)--(3,4);

     \draw[line width=1pt]
     (0,1)--(1,1)--(1,2)--(2,2)--(2,3)--(1,3)--(1,4);

     \draw[line width=1pt]
     (2,2)--(3,2)--(3,3)--(4,3);

     \draw[line width=1pt]
     (2,2)--(2,1)--(3,1)--(3,0);
    \end{tikzpicture}
    \end{subfigure}
    % third subfigure
    \begin{subfigure}{0.15\textwidth}
    \centering
     \begin{tikzpicture}[scale=0.4]
    % board
    \draw[thick] (0,0) rectangle (4,8);

    % Tile 1 (top left)
    \draw[line width=1pt]
    (0,7)--(1,7)--(1,6)--(2,6)--(2,7)--(3,7)--(3,8);

    \draw[line width=1pt]
    (2,6)--(3,6)--(3,5)--(4,5);

    \draw[line width=1pt]
    (0,4)--(1,4)--(1,5)--(2,5)--(2,6);

    \draw[line width=1pt]
    (1,4)--(4,4);
    %%%%
    \draw[line width=1pt]
    (0,3)--(1,3)--(1,2)--(2,2)--(2,3)--(3,3)--(3,4);

    \draw[line width=1pt]
    (2,2)--(3,2)--(3,1)--(4,1);

    \draw[line width=1pt]
    (0,0)--(1,0)--(1,1)--(2,1)--(2,2);

    \draw[line width=1pt]
    (1,4)--(4,4);

    \end{tikzpicture}
    \end{subfigure}
    % fourth subfigure
    \begin{subfigure}{0.15\textwidth}
    \centering
     \begin{tikzpicture}[scale=0.4]
    % board
    \draw[thick] (0,0) rectangle (4,8);

    % Tile 1 (top left)
    \draw[line width=1pt]
    (0,7)--(1,7)--(1,6)--(2,6)--(2,7)--(3,7)--(3,8);

    \draw[line width=1pt]
    (2,6)--(3,6)--(3,5)--(4,5);

    \draw[line width=1pt]
    (0,4)--(1,4)--(1,5)--(2,5)--(2,6);

    \draw[line width=1pt]
    (1,4)--(4,4);
    %%%%
    \draw[line width=1pt]
     (0,4)--(3,4);

     \draw[line width=1pt]
     (0,1)--(1,1)--(1,2)--(2,2)--(2,3)--(1,3)--(1,4);

     \draw[line width=1pt]
     (2,2)--(3,2)--(3,3)--(4,3);

     \draw[line width=1pt]
     (2,2)--(2,1)--(3,1)--(3,0);
    \end{tikzpicture}
    \end{subfigure}
    % fifth subfigure
    \begin{subfigure}{0.15\textwidth}
    \centering
     \begin{tikzpicture}[scale=0.4]
    % board
    \draw[thick] (0,0) rectangle (4,8);
    \draw[line width=1pt]
     (0,5)--(1,5)--(1,6)--(2,6)--(2,7)--(1,7)--(1,8);

    \draw[line width=1pt]
     (2,6)--(3,6)--(3,7)--(4,7);

    \draw[line width=1pt]
     (2,6)--(2,5)--(3,5)--(3,4);
    % Tiles in the second half
    \draw[line width=1pt]
    (0,3)--(1,3)--(1,2)--(2,2)--(2,3)--(3,3)--(3,4);

    \draw[line width=1pt]
    (2,2)--(3,2)--(3,1)--(4,1);

    \draw[line width=1pt]
    (0,0)--(1,0)--(1,1)--(2,1)--(2,2);
    %%%%%%
    \draw[line width=1pt]
    (0,4)--(4,4);
    \end{tikzpicture}
    \end{subfigure}
    % fifth subfigure
    \begin{subfigure}{0.15\textwidth}
    \centering
     \begin{tikzpicture}[scale=0.4]
    % board
    \draw[thick] (0,0) rectangle (4,8);

     \draw[line width=1pt] (0,3)--(1,3)--(1,2)--(2,2)--(2,1)--(1,1)--(1,0);
     \draw[line width=1pt] (2,2)--(3,2)--(3,1)--(4,1);
     \draw[line width=1pt] (0,4)--(1,4)--(1,5)--(2,5)--(2,2);

     \draw[line width=1pt] (2,3)--(3,3)--(3,4)--(4,4);
     \draw[line width=1pt] (0,7)--(1,7)--(1,6)--(2,6)--(2,5);

     \draw[line width=1pt] (2,6)--(3,6)--(3,5)--(4,5);

     \draw[line width=1pt] (2,6)--(2,7)--(3,7)--(3,8);
    \end{tikzpicture}
    \end{subfigure}
    \caption{Tilings for the $8\times 4$ board using $T$-tetrominoes.}
    \label{fig:8by4rot}
\end{figure}

The first few terms of the sequence for the generating function $F(8, T,t,4)$ are
\[ 1, 6, 84, 1182, 16644, 234390, 3300852, 46485102, 654638628, \ldots \]
In Figure~\ref{fig:8by4rot}, we display the $6$ tilings for the $8\times 4$ board for $n=1$.
The previous sequence does not currently appear in the OEIS. Nevertheless, it has been studied before. In 2004, Korn and Pak showed that the number of tilings of a $4k \times 4n$ board by $T$-tetrominoes can be expressed as an evaluation of the Tutte polynomial of the grid graph along a particular point~\cite{KornPak04}. In 2008, Merino evaluated these Tutte polynomials using the transfer-matrix method for fixed widths~\cite{Merino08}. Zeilberger's program provides an alternative computational approach that automatically produces the rational generating functions for each fixed width.

The next values of $k$ do not appear in the OEIS but are provided by Merino~\cite{Merino08}. Applying Zeilberger's algorithm to the rotations of the $T$-tetromino yields the sequences in Table~\ref{tab:t-tetro-seqs}. We include the generating functions and their corresponding recurrences in Appendix~\ref{appendix:t-tetro}.

\begin{table}[H]
    \centering
    \begin{tabular}{|c|c|c|}
     \hline 
       Board dimension  & First few terms of the sequence & OEIS \\\hline 
        $4\times 4n$ & $1,2, 6, 18, 54, 162, 486, 1458, 4374, 13122, 39366$ & \seqnum{A008776} \\ \hline 
        $8\times 4n$ & $1, 6, 84, 1182, 16644, 234390, 3300852, 46485102$ &  New \\ \hline 
        $12\times 4n$ & $1, 18, 1182, 78696, 5253822, 350950482, 23446010520$ &  New  \\ \hline 
        $16\times 4n$ & $1, 54, 16644, 5253822, 1668091536, 530454033510$ & New \\ \hline 
    \end{tabular}
    \caption{Sequences of tiling boards with rotations of the $T$-tetromino. The last three sequences do not appear in the OEIS but are studied by Merino~\cite{Merino08}. }
    \label{tab:t-tetro-seqs}
\end{table}

\subsection{Tiling boards with $L$-tetrominoes}\label{subsec:Ltetro}
\begin{figure}[H]
\begin{tikzpicture}[x=0.6cm,y=0.6cm]
%L-shape:
\redcell{0}{2}
\redcell{0}{1}
\redcell{0}{0}\redcell{1}{0}
%L-shape 2:
\redcell{3}{1}\redcell{4}{1}\redcell{5}{1}
\redcell{3}{0}
%L-shape 3:
\redcell{7}{2}\redcell{8}{2}
\redcell{8}{1}
\redcell{8}{0}
%L-shape 3:
\redcell{12}{1}
\redcell{10}{0}\redcell{11}{0}\redcell{12}{0}
\end{tikzpicture}
\caption{The set of all rotations of the $L$-tetromino.}
\label{fig:Ltetro-rot}
\end{figure}
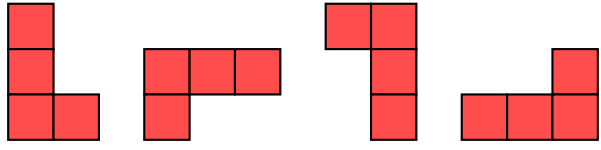

We now turn to the $L$-tetromino (see Figure~\ref{fig:Ltetro-rot}), which motivated the present work. Golomb and Klarner proved that a $k\times n$ board can be tiled by $L$-tetrominoes if and only if $k,n\geq 2$ and $kn$ is divisible by $8$~\cite{GolombKlarner}. Hence, for each $k\in \{5,6,7,8,9\}$, we consider the corresponding sequence $F(k,T,t,c)$ determined by the convention in Remark~\ref{rem:convention}.

As a first step, we show that Zeilberger's program automatically recovers the generating function recently obtained by B\v{e}lohoubek and Slav\'ik for $k=4$.
In~\cite[Theorem 1]{Ltetromino}, B\v{e}lohoubek and Slav\'ik prove that if $a_n$ denotes the number of tilings of a $2n\times 4$ board by $L$-tetrominoes, then \[a_n=a_{n-1}+3a_{n-2}-a_{n-3}-a_{n-4}, \quad n\geq 4,\]
and its generating function equals 
\begin{equation}\label{eq:genfunLtetro}
\frac{1-z^2}{z^4+z^3-3z^2-z+1}.
\end{equation}
Applying Zeilberger's program with $k=4$ in Maple and the set of tiles $T$ given by
\begin{align*}
T=\{ &\{ (0, 0), (0, 1), (0, 2), (1, 0)\}, \{(0, 0), (0, 1), (1, 1), (2, 1)\},\\
&\{(0, 0), (1, 0), (2, 0), (2, 1)\}, \{(0, 2), (1, 0), (1, 1), (1, 2)\}\},
\end{align*}
yields
\[\GFt(4,T,t)=\frac{1-t^{4}}{t^{8}+t^{6}-3t^{4}-t^{2}+1},\]
which can be rewritten as
\[ F(4,T,t,2)=\frac{1-t^{2}}{t^{4}+t^{3}-3t^{2}-t+1}.\]
By using the generating function $F(4,T,t,2)$, we obtain the first $12$ terms of the sequence, indexed from $0$ to $11$,
\[1, 1, 3, 5, 12, 23, 51, 103, 221, 456, 965, 2009.\]
This recovers Eq.~(\ref{eq:genfunLtetro}) of B\v{e}lohoubek and Slav\'ik and it corresponds to precisely the OEIS sequence \seqnum{A131322}, appearing in~\cite[Table 1]{Ltetromino}.
 
Having recovered the known generating function for the $2n\times 4$ case, we now apply the same computational framework to larger board widths. As an example, we give the generating function for the sequence corresponding to $5\times 8n$ boards.

\begin{proposition}\label{prop:5by8n-Ltetro}
Let $a_n$ denote the number of tilings of a $5\times 8n$ board by $L$-tetrominoes. The generating function is
\[ \sum_{n\geq 0}a_nt^n= \frac{2 t^{4}+2 t^{3}+t^{2}-1}{4 t^{4}+4 t^{3}+4 t^{2}+2 t -1}.\]
Consequently,
\[a_n=2a_{n-1}+4a_{n-2}+4a_{n-3}+4a_{n-4}, \quad n\geq 5,\]
with initial values $a_0=1, a_1=2, a_2=7, a_3=24$, and $a_4=86$.
The first terms of the sequence are
\[1,2,7,24,86,304,1076,3808, \ldots \]
\end{proposition}

Proposition~\ref{prop:5by8n-Ltetro} demonstrates how Zeilberger's program can be used to obtain additional tiling sequences. In Table~\ref{tab:Ltetro-seqs}, we summarize the sequences obtained from Zeilberger's program, several of which appear to be new to the OEIS. For comparison, we include the $k=4$ case previously studied. Observe that for $k=8$, there are lengths for which no tiling exists and the terms of the sequence are not increasing. The corresponding generating functions for each case are listed in Appendix~\ref{appendix:l-tetro}.
The same computational framework applies equally well to set of tiles that include reflections or to other free polyominoes. We expect many additional sequences arising from these classes, many of which are likely to be new to the OEIS.
We leave it as an open problem to determine bijections relating these sequences to existing combinatorial objects.

\begin{table}[H]
    \centering
    \begin{tabular}{|c|c|c|}
     \hline 
       Board dimension  & First few terms of the sequence & OEIS \\\hline 
        $4\times 2n$ & $1,1,3,5,12,23,51,103,221,456,965$ & \seqnum{A131322} \\ \hline 
        $5\times 8n$ & $1, 2, 7, 24, 86, 304, 1076, 3808, 13480$ &  New \\ \hline 
        $6\times 4n$ & $1,5,34,249,1882,14492,113029$ & New\\ \hline 
        $7\times 8n$ & $1, 14, 511, 18826, 699540, 26023308$ & New \\ \hline 
        $8\times n$ & $1,0,1,0,12,2,34,14,250,114, 956$ & New \\ \hline
        $9\times 8n$ & $1,114,42305,17416144,7258193991$ & New \\ \hline 
    \end{tabular}
    \caption{Sequences of tiling boards with rotations of the $L$-tetromino. For $k=8$, we include lengths for which no tiling exists.}
    \label{tab:Ltetro-seqs}
\end{table}

\section{Acknowledgments}
The author thanks her advisor, Dr. Doron Zeilberger, for feedback on earlier drafts. The author was supported by the NSF Graduate Research Fellowship Program under Grant No.~2233066 and in part by a Joel Lebowitz Summer Research Fellowship.

\bibliographystyle{plainurl}
\bibliography{references.bib}

\appendix
\section{Generating Functions for tiling boards}\label{appendix:generatingfuncs}
The generating functions in this section can also be found in the GitHub repository \href{https://github.com/marti310/TILINGSapps}{https://github.com/marti310/TILINGSapps}.

\subsection{Tiling with $T$-tetrominoes}\label{appendix:t-tetro}
The text file accompanying this section can be found in the GitHub repository~\cite{github} and is called \texttt{PaperTtetromino.txt}. The text file contains the generating functions, their corresponding sequences, and asymptotics that appear in this section.

\begin{proposition}
Let $a_n$ denote the number of tilings of a $8\times 4n$ board by $T$-tetrominoes. The generating function is
\[ \sum_{n\geq 0}a_nt^n= \frac{15t^2-10t+1}{27t^2-16t+1}.\]
Consequently,
\[a_n=16a_{n-1}-27a_{n-2}, \quad n\geq 3,\]
with initial values $a_0=1, a_1=6$, and $a_2=84$.
The first terms of the sequence are
\[1,6,84,1182,16644,234390,3300852,46485102,\ldots \]
\end{proposition}

\begin{proposition}
Let $a_n$ denote the number of tilings of a $12\times 4n$ board by $T$-tetrominoes. The generating function is
\[\sum_{n\geq 0}a_nt^n= \frac{4617t^4-4401t^3+1056t^2-69t+1}{6561t^4-6183t^3+1440t^2-87t+1}.\]
Consequently,
\[a_n=87a_{n-1}-1440a_{n-2}+6183a_{n-3}-
6561a_{n-4}, \quad n\geq 5,\]
with initial values $a_0=1, a_1=18, a_2=1182, a_3=78696$, and $a_4=5253822$.
The first terms of the sequence are
\[ 1, 18, 1182, 78696, 5253822, 350950482, 23446010520, \ldots \]
\end{proposition}

\begin{proposition}\label{prop:16by4n-ttetro}
Let $a_n$ denote the number of tilings of a $16\times 4n$ board by $T$-tetrominoes.
The generating function is $\sum_{n\geq 0}a_nt^n=A(t)/B(t)$, where
\begin{align*}
A(t)&= 679922958195 t^{10}-1140307639290 t^{9}+711589401621 t^{8}-216477728064 t^{7}\\
&\qquad +35414077878 t^{6}-3245223852 t^{5}+168003018 t^{4}-4796864 t^{3}+71031 t^{2}\\
&\qquad -474 t +1 ,
\end{align*}
and
\begin{align*}
B(t)&=847288609443 t^{10}-1419508671696 t^{9}+884318355441 t^{8}-268307318592 t^{7}\\
&\qquad +43712072010 t^{6}-3980340000 t^{5}+204075342 t^{4}-5739200 t^{3}+82899 t^{2}\\
&\qquad -528 t +1 .
\end{align*}
The first terms of the sequence are
\[ 1,54,16644,5253822,1668091536,530454033510, \ldots \]
\end{proposition}

\subsection{Tiling with $L$-tetrominoes}\label{appendix:l-tetro}

The text file accompanying this section can be found in the GitHub repository~\cite{github} and is called \texttt{PaperLpolyomino.txt}. The text file contains the generating functions, their corresponding sequences, and asymptotics that appear in this section. In addition, the text file contains additional sequences for tiling with $L$-pentominoes.

\begin{proposition}\label{prop:6by4n-Ltetro}
Let $a_n$ denote the number of tilings of a $6\times 4n$ board by $L$-tetrominoes. The generating function is $-A(t)/B(t)$, where
\[ 
A(t)=(t -1)^{2} (t^{9}-17 t^{8}+114 t^{7}-386 t^{6}+714 t^{5}-732 t^{4}+412 t^{3}-123 t^{2}+18 t -1),
\]
and
\begin{align*}
B(t)&=t^{12}-24 t^{11}+239 t^{10}-1280 t^{9}+4051 t^{8}-7944 t^{7}+9918 t^{6}-8001 t^{5}\\
&\qquad +4138 t^{4}-1330 t^{3}+251 t^{2}-25 t +1.
\end{align*}
The first terms of the sequence are
\[1, 5, 34, 249, 1882, 14492, 113029, 889719, 7050677, \ldots \]
\end{proposition}

\begin{proposition}\label{prop:7by8n-Ltetro}
Let $a_n$ denote the number of tilings of a $7\times 8n$ board by $L$-tetrominoes. The generating function is $-A(t)/B(t)$, where
\begin{align*}
A(t) &= 54 t^{25}-2805 t^{24}+61114 t^{23}-437586 t^{22}+756220 t^{21}-1683407 t^{20}+201818 t^{19}\\
&\qquad -4732871 t^{18}+5812686 t^{17}+6513779 t^{16}-5791118 t^{15}-15036852 t^{14}\\
&\qquad-14634004 t^{13} -854802 t^{12}+887014 t^{11}-5043234 t^{10}-3069842 t^{9}+558740 t^{8}\\
&\qquad +642802 t^{7}+30014 t^{6}-18954 t^{5}-2525 t^{4}-486 t^{3}-168 t^{2}-16 t +1 ,
\end{align*}
and
\begin{align*}
B(t) &= 3 t^{26}-572 t^{25}+20037 t^{24}-299100 t^{23}+1278161 t^{22}-3054798 t^{21}+5765025 t^{20}\\
&\qquad -3431380 t^{19} +23400340 t^{18}-15270754 t^{17}-4408598 t^{16}-17388698 t^{15}\\
&\qquad +20834748 t^{14}+39814896 t^{13}+28440798 t^{12}+16183108 t^{11}+10581385 t^{10}\\
&\qquad +6300318 t^{9}-880287 t^{8}-1336202 t^{7}-197778 t^{6}-1116 t^{5}-48 t^{4}\\
&\qquad +356 t^{3}+259 t^{2}+30 t -1 .
\end{align*}
The first terms of the sequence are
\[1, 14, 511, 18826, 699540, 26023308, 968374240, 36035363996, \ldots \]
\end{proposition}

\begin{proposition}\label{prop:8byn-Ltetro}
Let $a_n$ denote the number of tilings of a $8\times n$ board by $L$-tetrominoes. The generating function is $-A(t)/B(t)$, where
\begin{align*}
A(t) &=  4 t^{58}-5 t^{56}-69 t^{54}+34 t^{53}+36 t^{52}-28 t^{51}+334 t^{50}-194 t^{49}-181 t^{48}+196 t^{47}\\
&\qquad -560 t^{46} +356 t^{45}+490 t^{44}-396 t^{43}-22 t^{42}-114 t^{41}-240 t^{40}-44 t^{39}+1056 t^{38} \\
&\qquad -882 t^{37}-2511 t^{36} +1742 t^{35}-503 t^{34}+2046 t^{33}+9154 t^{32}-3982 t^{31}-2114 t^{30}\\
&\qquad -1838 t^{29}-16418 t^{28}+4992 t^{27} +4087 t^{26}+592 t^{25}+18004 t^{24}-3942 t^{23}-3377 t^{22}\\
&\qquad +80 t^{21}-12769 t^{20}+1978 t^{19}+1476 t^{18} -90 t^{17}+5910 t^{16}-602 t^{15}-344 t^{14}\\
&\qquad +20 t^{13} -1742 t^{12}+98 t^{11}+40 t^{10}-2 t^{9}+306 t^{8}-6 t^{7} -2 t^{6}-28 t^{4}+1 ,
\end{align*}
and
\begin{align*}
B(t) &= 8 t^{62}-13 t^{60}-12 t^{59}-133 t^{58}+96 t^{57}+103 t^{56}+138 t^{55}+594 t^{54}-770 t^{53}-406 t^{52}\\
& \qquad -730 t^{51} -758 t^{50} +2786 t^{49}+820 t^{48}+2328 t^{47}-672 t^{46}-6010 t^{45}+118 t^{44}\\ 
&\qquad -5192 t^{43}+1900 t^{42} +7148 t^{41}-7274 t^{40}+9842 t^{39}+2695 t^{38}-3132 t^{37}\\
&\qquad +28144 t^{36}-16280 t^{35}-14647 t^{34} -854 t^{33}-59822 t^{32}+21412 t^{31}+24741 t^{30}\\
&\qquad -464 t^{29}+79992 t^{28}-20802 t^{27}-22589 t^{26} +2654 t^{25}-70203 t^{24}+14230 t^{23}\\
&\qquad +11727 t^{22}-1922 t^{21}+40915 t^{20}-6502 t^{19}-3185 t^{18} +498 t^{17}-15672 t^{16}\\
&\qquad +1834 t^{15}+278 t^{14}-10 t^{13}+3807 t^{12}-278 t^{11}+59 t^{10}-16 t^{9} -543 t^{8}+16 t^{7}\\
&\qquad -15 t^{6}+2 t^{5}+39 t^{4}+t^{2}-1
.
\end{align*}
The first terms of the sequence are
\[1, 0, 1, 0, 12, 2, 34, 14, 250, 114, 956, 694, 5864, 4428, 25805, 25554, \ldots \]
\end{proposition}

We omit the generating function for the following result due to its size, but it can be found in the text file \texttt{PaperLpolyomino.txt} of the GitHub repository~\cite{github}.

\begin{proposition}\label{prop:9by8n-Ltetro}
Let $a_n$ denote the number of tilings of a $9\times 8n$ board by $L$-tetrominoes. The generating function is $-A(t)/B(t)$ where both $A(t)$ and $B(t)$ are polynomials of degree $146$.
The first terms of the sequence are
\[1, 114, 42305, 17416144, 7258193991, 3032714594518, 1267783654718834, \ldots \]
\end{proposition}

\end{document}